\documentclass{gtart}

\def\ifplaintex{\expandafter\ifx\csname documentclass\endcsname\relax}

\ifplaintex 
\else
\fi

\expandafter\ifx\csname beginpicture\endcsname\relax
\expandafter\ifx\csname documentclass\endcsname\relax
\input pictex \else
\input prepictex \input pictex \input postpictex \fi\fi

\def\gt{{\mathsurround=0pt\it $\cal G\mskip-2mu$eometry \&\ 
$\cal T\!\!$opology}}        %  journal title in recommended style

\def\gtp{{\mathsurround=0pt\it $\cal G\mskip-2mu$eometry \&\ 
$\cal T\!\!$opology $\cal P\!$ublications}}  % GT publications

\def\volumenumber#1{\def\thevolumenumber{#1}}
\def\papernumber#1{\def\thepapernumber{#1}}
\def\volumeyear#1{\def\thevolumeyear{#1}}

\def\pagenumbers#1#2{\def\startpage{#1}\def\finishpage{#2}}
\def\published#1{\def\publishdate{#1}}
\def\proposed#1{\def\theproposer{#1}}
\def\seconded#1{\def\theseconders{#1}}
\def\received#1{\def\receiveddate{#1}}
\def\revised#1{\def\reviseddate{#1}}
\def\accepted#1{\def\accepteddate{#1}}

\long\def\asciiabstract#1{\long\def\theasciiabstract{#1}}

\def\shorttitle#1{\def\theshorttitle{#1}}

\let\\\par\let\thevolumenumber\relax\let\thepapernumber\relax
\let\thevolumeyear\relax\let\thesamplenumber\relax\let\startpage\relax
\let\finishpage\relax\let\publishdate\relax\let\receiveddate\relax
\let\reviseddate\relax\let\accepteddate\relax\let\theasciititle\relax
\let\theasciiauthors\relax
\let\theasciiabstract\relax
\let\theasciiemail\relax\let\theshortauthors\relax\let\theshorttitle\relax

\long\def\maketitlep{   % start of definition of \maketitlep

\count0=\startpage

\gt\hfill      %   Journal title (top left) 
\beginpicture
\setcoordinatesystem units <0.33truein, 0.33truein> point at 2.2 0.9
\setplotsymbol ({$\cal G$})
\plotsymbolspacing=9truept
\circulararc 315 degrees from 0 1 center at 0 0
\setplotsymbol ({$\cal T$})
\circulararc 315 degrees from 1 -1 center at 1 0
\endpicture
\break
{\small\ifx\thesamplenumber\relax % sample?  
Volume \else Sample
\fi\thevolumenumber\ (\thevolumeyear)
\startpage--\finishpage\nl
Published: \publishdate}
\vglue 0.5truein plus 0.4fil minus 0.1truein

{\parskip=0pt\leftskip 0pt plus 1fil\def\\{\par\smallskip}{\ifplaintex\large
\else\Large\fi\bf\thetitle}\par\medskip}   

\vglue 0pt plus 0.1fil 

{\parskip=0pt\leftskip 0pt plus 1fil\def\\{\par}{\sc\theauthors}
\par\medskip}

\vglue 0pt plus 0.1fil 

{\small\parskip=0pt\let\newline\\
{\leftskip 0pt plus 1fil\def\\{\par}{\sl\theaddress}\par}
\expandafter\ifx\theemail\relax    % email address?
\relax\else\vglue 5pt plus 0.02fil minus 2pt\def\\{\stdspace{\rm 
and}\stdspace} 
\cl{Email:\stdspace\tt\theemail}\fi
\ifx\theurl\relax                  % URL given?
\relax\else\vglue 5pt plus 0.02fil minus 2pt\def\\{\stdspace{\rm 
and}\stdspace}
\cl{URL:\stdspace\tt\theurl}\fi\par}

\vglue 7pt plus 0.3fil minus 3pt

{\bf Abstract}
\vglue 5pt plus 0.1fil minus 2pt

\theabstract

\vglue 7pt plus 0.3fil minus 3pt

{\bf AMS Classification numbers}\quad Primary:\quad \theprimaryclass

Secondary:\quad \thesecondaryclass

\vglue 5pt plus 0.3fil minus 2pt

{\bf Keywords:}\quad \thekeywords

\vglue 10pt plus 0.5fil minus 5pt

{\small  Proposed: \theproposer\hfill Received: \receiveddate\nl
Seconded: \theseconders\hfill 
\ifx\reviseddate\relax                         % paper revised?
Accepted: \accepteddate                        % no
\else
Revised: \reviseddate                          % yes
\fi}
\eject
}       %  end of definition of \maketitlep

\let\maketitlepage\maketitlep
\let\maketitle\maketitlepage

\font\phead=cmsl9 scaled 950
\font=cmsl9 scaled 1050
\font\pnum=cmbx10 scaled 913
\font\lnum=cmbx10 
\font\pfoot=cmsl9 scaled 950
\font=cmsl9 scaled 1050
\ifplaintex
\headline{\vbox to 0pt{\vskip -4.5mm\line{\small\phead\ifnum
\count0=\startpage ISSN 1364-0380 (on line)
1465-3060 (printed) \hfill {\pnum\folio}\else\ifodd\count0\def\\{ }% 
\ifx\theshorttitle\relax\thetitle\else\theshorttitle\fi\hfill{\pnum\folio}
\else\def\\{ and }{\pnum\folio}\hfill\ifx\theshortauthors\relax\theauthors
\else\theshortauthors\fi\fi\fi}\vss}}
\footline{\vbox to 0pt{\vglue 0mm\line{\small\pfoot\ifnum\count0=\startpage
\copyright\ \gtp\hfill\else
\gt, Volume \thevolumenumber\ (\thevolumeyear)\hfill\fi}\vss
}}
\else
\makeatletter
\def\@oddhead{{\small\lhead\ifnum\count0=\startpage ISSN 1364-0380 (on line)
1465-3060 (printed) \hfill {\lnum\number\count0}\else\ifodd\count0
\def\\{ }\ifx\theshorttitle\relax \thetitle \else\theshorttitle\fi\hfill
{\lnum\number\count0}\else\def\\{ and }{\lnum\number\count0}
\hfill\ifx\theshortauthors\relax 
\theauthors\else\theshortauthors\fi\fi\fi}}\def\@evenhead{@oddhead}
\def\@oddfoot{\small\lfoot\ifnum\count0=\startpage\copyright\ \gtp\hfill\else
\gt, Volume \thevolumenumber\ (\thevolumeyear)\hfill\fi}
\def\@evenfoot{@oddfoot}
\makeatother
\fi

\newwrite\gtoutfile
\long\gdef\makeheadfile{  %%% start of definition of \makeheadfile
{\def\\{, }\def\s{ }
\immediate\openout\gtoutfile head.xxx
\immediate\write\gtoutfile{To: math@arxiv.org}
\immediate\write\gtoutfile{Subject: put}
\immediate\write\gtoutfile{--text follows this line--}
\immediate\write\gtoutfile{Proxy-for: \ifx\theasciiauthors\relax
\theauthors\else\theasciiauthors\fi\s<\ifx\theasciiemail\relax\theemail\else\theasciiemail\fi>}
\immediate\write\gtoutfile{\noexpand\\}
\immediate\write\gtoutfile{Authors: \ifx\theasciiauthors\relax
\theauthors\else\theasciiauthors\fi}
{\def\\{ }\immediate\write\gtoutfile{Title: \ifx\theasciititle\relax
\thetitle\else\theasciititle\fi}}
\immediate\write\gtoutfile{Subj-class: GT}
\immediate\write\gtoutfile{MSC-class: \theprimaryclass\ifx\thesecondaryclass\relax\else, \thesecondaryclass\fi}
\immediate\write\gtoutfile{Journal-ref: Geom. Topol. \thevolumenumber\s
(\thevolumeyear) \startpage-\finishpage}
\immediate\write\gtoutfile{Comments: Published in Geometry and Topology at}
\immediate\write\gtoutfile{    http://www.maths.warwick.ac.uk/gt/GTVol\thevolumenumber/paper\thepapernumber.abs.html}
\immediate\write\gtoutfile{\noexpand\\}
\immediate\write\gtoutfile{}
\ifx\theasciiabstract\relax
\immediate\write\gtoutfile{\theabstract}\else
\immediate\write\gtoutfile{\theasciiabstract}\fi
\immediate\write\gtoutfile{}
\immediate\write\gtoutfile{\noexpand\\}
\immediate\write\gtoutfile{}
\immediate\write\gtoutfile{<uuencoded .tar.gz file here>}
\immediate\write\gtoutfile{}
\immediate\closeout\gtoutfile}}  %%% end of definition of \makeheadfile

\def\maketitlepage{\maketitlep\makeheadfile}
\let\maketitle\maketitlepage

\volumenumber{4}\papernumber{10}\volumeyear{2000}
\pagenumbers{293}{307}
\proposed{Joan Birman}
\seconded{Shigeyuki Morita, Walter Neumann}
\received{24 November 1999}
\revised{28 September 2000}
\accepted{3 August 2000}
\published{10 October 2000}

\usepackage{amsmath,amsfonts}
\usepackage{graphicx,labelfig}

\newtheorem{thm}{Theorem}[section]
\newtheorem{lemma}[thm]{Lemma}
\newtheorem{sublemma}{Claim}
\newtheorem{defn}{Definition}[section]

\newenvironment{oneshot}[1]%
{\begin{trivlist}\item[\hskip\labelsep{\bf #1}]\sl}%
{\end{trivlist}}

\newcommand{\pf}{\smallskip \noindent \textbf{Proof}\qua}

\newcommand{\pffive}{\smallskip \noindent \textbf{Proof of Theorem 5.1}\qua}

\newcommand{\rmk}{\noindent \textbf{Remark}\qua}

\newcommand{\F}{F_p(X,n)}

\newcommand{\By}{F(Y,n)}
\newcommand{\Bdy}{F^\partial(Y,n)}
\newcommand{\Fy}{F_p(Y,n)}
\newcommand{\Fdy}{F_p^\partial(Y,n)}

\newcommand{\M}{M(\Sigma_2)}
\newcommand{\Ms}{M(S^2,6)}

\newcommand{\Mpxn}{M_p(X,n)}
\newcommand{\Myn}{M(Y,n)}
\newcommand{\Mpyn}{M_p(Y,n)}
\newcommand{\My}{M(Y)}
\newcommand{\Myd}{M^\partial(Y)}

\newcommand{\Mynd}{M^\partial(Y,n)}
\newcommand{\Mpynd}{M^\partial_p(Y,n)}

\newcommand{\MSnp}{M_p(S^2,n)}

\newcommand{\Mplus}{M^\partial(\Dp,1)}
\newcommand{\Mminus}{M_p^\partial(\Dm,k)}
\newcommand{\MSkp}{M_p(S^2,k+1)}

\newcommand{\Msfourp}{M_p(S^2,4)}

\newcommand{\Mccv}{M({\Ttwo})} %replacing this also
\newcommand{\Mcc}{M^\partial({\Ttwo})}  %replacing with \Mcc

\newcommand{\Mtp}{M_p(T,2)}
\newcommand{\Mt}{M(T,1)}

\newcommand{\Mtorone}{M^\partial(T_1 - 1)}
\newcommand{\Mtortwo}{M^\partial(T_2 - 1)}

\newcommand{\Br}{Br(S^2,n)}
\newcommand{\Brus}{Br(S^2,6)}
\newcommand{\Brt}{Br(\TT)}

\newcommand{\pT}{\pi_1(T - p_1)}

\newcommand{\slthree}{Sp(4,\Z_3)}

\newcommand{\Z}{\mathbb{Z}}

\newcommand{\ypoints}{\{y_1,\ldots,y_n\}}

\newcommand{\Xd}{\overline{X-D}}
\newcommand{\Xdd}{\overline{X-D'}}
\newcommand{\Dd}{\overline{D'-D}}

\newcommand{\Xnm}{X\setminus{\{x_1,\ldots, x_{n-1}\}}}
\newcommand{\xptnm}{x_1,\ldots,x_{n-1}}

\newcommand{\Ttwo}{T-2}
\newcommand{\torone}{T_1 - 1}
\newcommand{\tortwo}{T_2 - 1}

\newcommand{\fs}{F^s}
\newcommand{\fu}{F^u}

\newcommand{\TT}{\Sigma_2}

\newcommand{\Dp}{D^{+}}
\newcommand{\Dm}{D^{-}}

\newcommand{\TDm}{T_{\partial \Dm}}
\newcommand{\TDp}{T_{\partial \Dp}}

\newcommand{\tw}{\langle \Tbp,\Tbm \rangle}

\newcommand{\Tbp}{T_{\partial1}}
\newcommand{\Tbm}{T_{\partial2}}

\begin{document}

\title{Normal all pseudo-Anosov subgroups of\\mapping class groups}
\shorttitle{All pseudo-Anosov subgroups}
\author{Kim Whittlesey}

\address{Department of Mathematics, The Ohio State University\\231
W 18th Avenue, Columbus, OH 43210, USA}
\email{whittle@math.ohio-state.edu}

\begin{abstract}                      
  
  We construct the first known examples of nontrivial, normal, all
  pseudo-Anosov subgroups of mapping class groups of surfaces.
  Specifically, we construct such subgroups for the closed genus two
  surface and for the sphere with five or more punctures. Using the
  branched covering of the genus two surface over the sphere and
  results of Birman and Hilden, we prove that a reducible mapping
  class of the genus two surface projects to a reducible mapping class
  on the sphere with six punctures. The construction introduces
  ``Brunnian'' mapping classes of the sphere, which are analogous to
  Brunnian links.

\end{abstract}

\asciiabstract{We construct the first known examples of nontrivial,
  normal, all pseudo-Anosov subgroups of mapping class groups of
  surfaces.  Specifically, we construct such subgroups for the closed
  genus two surface and for the sphere with five or more
  punctures. Using the branched covering of the genus two surface over
  the sphere and results of Birman and Hilden, we prove that a
  reducible mapping class of the genus two surface projects to a
  reducible mapping class on the sphere with six punctures. The
  construction introduces "Brunnian" mapping classes of the sphere,
  which are analogous to Brunnian links.}

\primaryclass{57M60}\secondaryclass{20F36, 57N05}

\keywords{Mapping class group, pseudo-Anosov, Brunnian}

\maketitlepage

\section
{Introduction} In 1985, R\,C Penner, D\,D Long, and J\,D McCarthy
asked the following question~\cite{kn:lon}: does there
exist a nontrivial, normal subgroup of the mapping class group of a
surface, all of whose nontrivial elements are pseudo-Anosov?  This
paper gives a positive answer for the case of a closed genus two
surface and for the case of the sphere with five or more punctures. The
approach is to construct a normal subgroup that avoids reducible and
periodic elements; by work of Thurston, the remaining elements must
all be pseudo-Anosov.  

A Brunnian link is defined to be a nontrivial link such that every
proper sublink is trivial. In section two we will define ``Brunnian''
mapping classes of the sphere and ``Brunnian'' sphere braids.  Roughly
speaking, a Brunnian sphere braid is analogous to a Brunnian link: it
cannot ``lose'' a strand without becoming trivial; similarly,
``Brunnian mapping classes of the n--punctured sphere'' are those that
become trivial if \textit{any} of the punctures is ``forgotten.''  The
subgroup consisting of Brunnian mapping classes is a normal subgroup
of the mapping class group of the n--punctured sphere.  In section
four we show that for $n \geq 5$, this subgroup has no reducible or
periodic elements, thus proving

\begin{thm} 
  For $n \geq 5$, The Brunnian subgroup $Br(S^2,n)$ is a
 nontrivial subgroup of $M(S^2,n)$ all of whose nontrivial elements are
 pseudo-Anosov.
\end{thm}

We use work of Birman and Hilden to lift the Brunnian subgroup of the
mapping class group of the 6--punctured sphere to a ``Brunnian
subgroup'' of the mapping class group of the closed genus two surface.
We prove in section five that if a genus two mapping class is
reducible, then the projection to the punctured sphere must also be
reducible.  The Brunnian subgroup for the genus two surface thus has
no reducible elements except for the involution that gives the
branched covering. To get rid of this element, we intersect with the
kernel of the usual $\Z_3$ homology representation.

\begin{thm}
The intersection of the Brunnian subgroup $Br(\TT)$ and the kernel of the
$\Z_3$ homology representation, $\rho\co  \M \to \slthree$, is a
nontrivial normal subgroup of $\M$, all of whose nontrivial elements
are pseudo-Anosov.
\end{thm}

\section{Background and definitions}

Let $Y$ be an orientable surface, possibly with boundary, and let
$\ypoints$ be an n--element subset of $Y$. We begin by introducing a
notation that will let us easily distinguish between several different
mapping class groups for $Y$.  

\begin{defn}\rm

\noindent We define
\begin{verse}
 
  $\By$ to be the group of all orientation preserving homeomorphisms
  ${h\co  Y \to Y}$ that induce the identity permutation on the set of
  boundary components of $Y$ and take the set of points $\ypoints$ to
  itself, and \\
  $\Myn   = \pi_0(\By)$ to be the corresponding mapping class group,\\
  \vspace{.5cm}
  
  $\Bdy$, with a superscript $\partial$, to be the group of all
  orientation preserving homeomorphisms ${h\co  Y \to Y}$ which fix the
  boundary pointwise and take the set of points $\ypoints$ to itself,
  and
  \\
  $\Mynd  = \pi_0(\Bdy$ to be the corresponding mapping class group,\\
  
  \vspace{.5cm} $\Fy$, with a subscript $p$ (in analogy with pure
  braids), to be the group of all orientation preserving
  homeomorphisms ${h\co  Y \to Y}$ that induce the identity permutation
  on the set of boundary components of $Y$, and also fix the $y_i$
  pointwise, ie, $h(y_i) = y_i$ for all $1 \leq i \leq n$,\\
  
  $\Mpyn  = \pi_0(\Fy )$ to be the corresponding mapping class group,\\
  
  \vspace{.5cm} $\Fdy$, with both a superscript $\partial$ and a
  subscript $p$, to be the group of all orientation preserving
  homeomorphisms ${h\co  Y \to Y}$ which fix the boundary pointwise and
  fix the $y_i$ pointwise, and\\
  
  $\Mpynd = \pi_0(\Fdy)$ to be the corresponding mapping class group.\\

\end{verse}
\end{defn}

Note that the elements of these mapping class groups \textit{are not
  permitted to permute the boundary curves.} A superscript $\partial$
denotes that the homeomorphisms and isotopies restrict to the identity
on the boundary curves. A subscript $p$ denotes that the
homeomorphisms do not permute the points $y_1, \ldots, y_n$.  If the
surface $Y$ has no boundary, then $\Mynd = \Myn$; if there are no
specified points, then we write simply $\My$ or $\Myd$.  The natural
homomorphism from $\Mynd$ to $\Myn$ has kernel generated by Dehn
twists along curves parallel to the boundary curves.

The mapping class group most commonly found in the literature would
have $M(Y)$ as a subgroup; in these larger mapping class groups,
mapping classes may permute the boundary components, and isotopies do
not restrict to the identity on the boundary components.

The notation $[f]$ is used to denote the mapping class represented by
a member $f$ of any of these homeomorphism groups.

\begin{thm}[Thurston~\cite{kn:flp}]
  If $Y$ is an orientable surface of negative Euler characteristic,
  then a mapping class $[h]$ of $\My$ must be either
\begin{enumerate}
\item[\rm(i)]{periodic, ie, $[h]^n = 1$ for some $n \neq 0$, or}
\item[\rm(ii)]{reducible, ie, having a representative that leaves
    invariant a non-empty collection of disjoint essential simple
    closed curves, none of which are isotopic to boundary curves, or }
\item[\rm(iii)]{pseudo-Anosov, ie, having a representative that leaves
 invariant a transverse pair of singular foliations $\fs, \fu$ with
 transverse measures $\mu^s, \mu^u$ such that $h(\fs, \mu^s) = (\fs,(
1/\lambda) \mu^s)$ and $h(\fu, \mu^u) = (\fu, \lambda\mu^u)$ for some
 $\lambda > 1$.}
\end{enumerate}
\end{thm}

In knot theory, a Brunnian link is defined to be a nontrivial link
such that every proper sublink is trivial~\cite[page 67]{kn:rolf}. In
\cite{kn:lev} Levinson introduces disk braids with an analogous
property: a ``decomposable braid'' is one such that if any single
arbitrary strand is removed, the remaining braid is trivial. We will
define a
\textit{Brunnian sphere braid} analogously: removal of
\textit{any} strand gives a trivial braid. See Figure 1. Note that
these braids must always be pure. Similarly, a \textit{Brunnian
mapping class of the sphere} will be one where if we ``forget'' any of
the points, there is an isotopy to the identity fixing all of the
other points.  Note that in order to ``remove'' a strand or ``forget''
a point in a well defined way, we need to assume that the permutation
induced by the mapping class actually fixes that strand or point.

\begin{figure}[ht!]
\centerline{\includegraphics[width=\hsize]{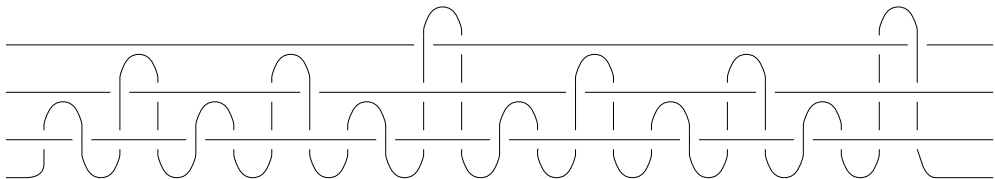}}
\caption{A  Brunnian sphere braid with four strands}
\end{figure}

Formally, then, let 
\begin{verse}
  $B(S^2,n,y_i)$ be the subgroup of
  $B(S^2,n)$ consisting of homeomorphisms such that $h(y_i) = y_i$, ie, the induced permutation of the points leaves $y_i$ fixed, and let\\
  $M(S^2,n,y_i) = \pi_0(B(S^2,n,y_i))$.\\
\end{verse}
Note that $M(S^2,n,y_i)$ is a subgroup of $M(S^2,n)$ of index $n$.
Since a homeomorphism in $B(S^2,n,y_i)$ leaves the set of $n-1$ points
$\{y_1,\dots,y_n\}\setminus y_i$ invariant, there is a well defined
inclusion
$$
B(S^2,n,y_i) \to B(S^2,n-1),
$$
which induces a homomorphism
$$
L'_i\co  M(S^2,n,y_i) \to M(S^2,n-1).
$$
We now consider the intersection of all $n$ of these kernels,
$\bigcap_{i=1}^{n} \ker{(L'_i)}$. This is clearly a normal subgroup of
$M(S^2,n)$. An element in this subgroup must act as the identity
permutation on the $n$ points.  

Alternatively, if we define $L_i$ to be the restriction of $L'_i$ to
$M_p(S^2,n)$,
$$
 L_i\co  M_p(S^2,n) \to M_p(S^2,n-1),
$$
 then
$$
\bigcap_{i=1}^{n} \ker{(L_i)} = \bigcap_{i=1}^{n} \ker{(L'_i)}.
$$

\begin{defn}\rm
We will call the maps $L_i$ the \textit{forget $y_i$ maps,}
and define the \textit{n--strand Brunnian subgroup of the sphere}, $\Br$, to be
$$
\Br = \bigcap_{i=1}^{n} \ker{(L_i)} = \bigcap_{i=1}^{n} \ker{(L'_i)}.
$$
Elements of this subgroup will be called \textit{Brunnian mapping classes} of
$M(S^2,n)$.
\end{defn}

We now define Brunnian mapping classes of the genus two torus. Let
$\TT$ be the closed orientable genus two surface, and let $\M$ be its
mapping class group. There is a 2--fold branched covering $\pi\co  \TT
\to S^2$; it has six branch points which we will denote by $y_1,\dots,
y_6 \in S^2$.  The group of covering transformations is generated by a
hyper-elliptic involution $i$, and $[i]$ is central in $\M$.  Let $\Ms$
be the mapping class group of the sphere, with the branch points as
the specified points.
\begin{lemma}{\rm(Birman and Hilden \cite[pages 183--189]{kn:birm})}\qua
 There is an exact sequence
$$
 1 \to \langle [i] \rangle \to \M \stackrel{p}{\to} \Ms \to 1\
$$
where $p$ takes the standard generating Dehn twists
$[\Delta_1],\dots,[\Delta_5]$ to the standard generators of $\Ms$,
$[\sigma_1], \dots, [\sigma_5]$, which switch adjacent branch points
on $S^2$.
\end{lemma}

\begin{defn}\rm
The \textit{Brunnian subgroup $Br(\TT)$ of $\M$} is the preimage of
the six strand Brunnian subgroup of the sphere:
$$ 
%R = 
Br(\TT) = p^{-1} (\Brus) 
%= p^{-1}( \bigcap_{i=1}^{6} \ker{(L_i)}).
$$
Elements of $Br(\TT)$ will be called \textit{Brunnian mapping classes of $\TT$.}
Note that $Br(\TT)$ is a normal subgroup of $\M$.
\end{defn}

Finally, we discuss homology with $\Z_3$ coefficients.  Since
$H_1(\TT, \Z_3) = (\Z_3)^4$, there is a representation
$$\rho\co  \M \to \slthree$$
taking a mapping class to its induced map on homology.

We now restate the main theorem of this paper.

\begin{oneshot} {Theorem 1.2}
The intersection of the Brunnian subgroup $Br(\TT)$ and the kernel of the
$\Z_3$ homology representation, $\rho\co  \M \to \slthree$, is a
nontrivial normal subgroup of $\M$, all of whose nontrivial elements
are pseudo-Anosov.
\end{oneshot}

\rmk The methods used in this paper will not work to produce normal all
pseudo-Anosov subgroups for higher genus surfaces, since those mapping
class groups contain no central element to act as a covering
transformation.  However, we conjecture that normal all pseudo-Anosov
subgroups do exist for higher genus surfaces.

\section{Preliminary lemmas} 
We need the following fairly basic lemma, whose proof we include for
the sake of completeness.

\begin{lemma}
  Given a closed orientable surface $X$ with specified points\break $x_1,
  \dots, x_n$, and a closed disk $D \subset X$ containing $x_n$ but
  none of $x_1,...,x_{n-1}$, the map $$M_p(\Xd,n-1) \to \Mpxn$$
  induced by inclusion is an isomorphism.
 \end{lemma}
 
 \pf 
Recall that $M_p(\Xd,n-1)$ is defined to be the group of isotopy
 classes of self homeomorphisms of $\Xd$ that fix each of $x_1,\ldots,
 x_{n-1}$.  Any homeomorphism or isotopy of $\Xd$ fixing $x_1,\ldots,
 x_{n-1}$ can be extended to one on $X$ fixing $x_1,\ldots, x_n$.  To
 prove that the above map is one-to-one, we must show: given a
 homeomorphism $f\co  X \to X$ that fixes each of $\xptnm$ and leaves $D$
 invariant, and an isotopy $h_t\co X \to X$ fixing $x_1,\ldots, x_n$ and
 going from the identity to $f$, there is an isotopy $H_t\co \Xd \to \Xd$
 going from the identity on $\Xd$ to the restriction of $f$ to $\Xd$.

 Choose a closed disk $D'$ such that $D \subset
\mathrm{int}D'$, $x_1, \dots, x_{n-1} \notin D'$, and further \linebreak[1]
${h_t(\Xdd) \cap D = \phi}$ for all $t.$ Let $g_t\co  D \cup \Xdd \to X$ be
an isotopy that is constant on $D$ and agrees with $h_t$ on $\Xdd$.
Using an isotopy extension theorem for topological
manifolds~\cite{kn:edkr}, $g_t$ can be extended to an isotopy
$G_t\co X\to X$.  $G_0$ is the identity map; thus by construction,
$G_1|_{\Xd}$ is trivial in $M_p(\Xd,n-1)$.

Unfortunately, $G_1|_{\Xd}$ may not agree with $f|_{\Xd}$. They are
the same on $\Xdd$, but $f|_{\Dd}$ may not equal $G_1|_{\Dd}$.
However, on $\Dd$, an annulus, \textit{any} two functions are isotopic
to each other by an isotopy which fixes $\partial D'$ pointwise but
can vary along $\partial D$.  If we attach such an isotopy between
$G_1|_{\Dd}$ and $f|_{\Dd}$ to the constant isotopy on $\Xdd$, we get
an isotopy $J_t\co  \Xd \to \Xd$ going from $G_1|_{\Xd}$ to
$f|_{\Xd}$. 
\endproof

\begin{lemma}{\rm(Birman \cite[page 217]{kn:birm2})}\qua
  Let $X$ be a a closed orientable surface with specified points
  $x_1,\ldots,x_n$, and let $L_n\co  \Mpxn \to M_p(X,n-1)$ be the
  ``forget $x_n$ map''. Then there is a long exact sequence ending
$$
\to \pi_1(F_p(X,n-1))  \to \pi_1(\Xnm) \stackrel{d}{\to} \Mpxn 
\stackrel{L_n}{\to} M_p(X,n-1) \to 1.
$$
Moreover, $\ker{d}$ is contained in the center of $\pi_1(\Xnm).$
\end{lemma}

\rmk
Since $\pi_1(\Xnm)$ has trivial center unless $X = S^2$ and 
$n=3$, or $X = T$ and  $n=1$, in all but these cases we get a short exact
sequence 
$$
 1 \to \pi_1(\Xnm) \to \Mpxn \stackrel{L_n}{\to} M_p(X,n-1) \to 1.
$$

\pf The proof follows the lines of Lemmas 4.1.1 and 4.2.1 in
[B].  Namely, first one must show that the evaluation map 
$$
 \epsilon\co  F_p(X,n-1) \to \Xnm
$$
taking $h$ to $h(x_n)$ is a locally trivial fibering map with fiber
$\F$.  The exact sequence above is then the exact homotopy sequence
of the fibering.\break  Finally, elements of $\ker{d}$ can be shown to
commute with everything in\break${\pi_1(\Xnm)}$.
\endproof

\begin{section}{Periodic and reducible Brunnian mapping classes of $S^2$ are 
trivial}

We first consider periodic maps in $M(S^2,n)$. Lemma 3.2 tells us that
for $n \geq 4$, there is an exact sequence
$$
1 \to \pi_1(S^2\setminus{\{x_1,...,\widehat{x_i},...,x_n\}}) \to \MSnp
\stackrel{L_i}{\to} M_p(S^2,n-1) \to 1
$$
for each $L_i$, the ``forget about $y_i$'' map. Therefore, when $n
\geq 4$, the subgroup of Brunnian sphere mapping classes is a free
group, and has {\em no nontrivial periodic elements.\/}
  
\begin{defn}\rm
  A reducible homeomorphism in $F(S^2,n)$ is one that leaves invariant
  a collection of disjoint essential closed curves, none of which is
  isotopic to a simple loop around a single point.  A reducible
  element of $M(S^2,n)$ is one with a reducible representative. 

\end{defn}

\begin{thm}
 For $n\geq 5$, an element of $Br(S^2,n)$ that is also reducible must
be trivial.
\end{thm}

\pf Let $[g]$ be a reducible element of $Br(S^2,n)$, for some $n \geq
5$, and let $g$ be a representative of $[g]$ that leaves invariant a
collection of disjoint essential closed curves in $S^2$.  Let $E$ be
one these curves.  $E$ separates the sphere into two open hemispheres,
each containing at least two of the $n$ points $y_1, \dots,y_n$.  Some
power $h = g^r$ leaves the curve $E$ itself invariant.  Since $[h]$ is
Brunnian, it does not permute the $n$ points, and so we can isotope
$h$ so that that it restricts to the identity on $E$. We will show
that $[h]$ is trivial. Then, since $Br(S^2,n)$ contains no nontrivial
periodic elements, $[g]$ must also be trivial.

Let $\Dp$ and $\Dm$ be the closures of the two components of $S^2
\setminus E$.  We will refer to these as the upper and lower
hemispheres, respectively.  If a homeomorphism $h$ fixes $E$ and the
$n$ points $y_1, \dots, y_n$, then it restricts to homeomorphisms of
$\Dp$ and $\Dm$.  For convenience, we define the following notation:
if $h|\Dp = f$ and $h|\Dm = g$, we write $h = (f,g)$.  Similarly,
given any two homeomorphisms $f$ and $g$ on the hemispheres that fix
the points and $E$, $(f,g)$ denotes the homeomorphism of the sphere
gotten by attaching $f$ and $g$ along $E$.

As we are interested in Brunnian mapping classes, let us consider what
happens when we ``forget'' all but one of the points on the top
hemisphere. Let $k$ be the number of the $y_i$ on the bottom
hemisphere, so that there are $n-k$ points on the top hemisphere.
Forgetting all but one of the $y_i$ from the top hemisphere will
leave a total of $k+1$ points to be fixed under isotopies. 

Consider a Brunnian mapping class $[h]$ with a representative $h$
that fixes $E$.  When we forget all but one point of the upper
hemisphere, the image of $[h]$ in $\MSkp$ must be trivial.  
Let $h_+ = h|\Dp$ and $h_- = h|\Dm$, so that $h = (h_+,h_-)$.
\begin{sublemma}
$ [(h_+, h_-)] = [(1,h_-)]$ as elements of $\MSkp$.
\end{sublemma}
\pf
$\Mplus = {1}$, so there is an isotopy of the disk $\Dp$, fixing the
one point and the boundary, taking $h_+$ to the identity.  Attach this
isotopy to the constant isotopy on $\Dm$.
\endproof
\begin{sublemma}
The homomorphism $j \co \Mminus \to \MSkp$ taking $[f] \to [(1,f)]$ has
kernel generated by a single Dehn twist $[\TDm]$ along a simple closed curve
parallel to the boundary.
\end{sublemma}
\pf
 $\MSkp \cong M_p(\Dm,k)$, by Lemma 3.1. However, the kernel of the
map $\Mminus \to M_p(\Dm,k)$ is generated by the Dehn twist along
a curve parallel to the boundary of $\Dm$.
\endproof

Since $[(1,h_-)] = [(h_+,h_-)] = [h] = 1$ as elements of $\MSkp$,
$[h_-]$ is a power of $[\TDm]$ in $\Mminus$, and so $h_-$ is isotopic
to a power of $\TDm$ by an isotopy fixing the boundary and all of the
original $y_i$ that were on $\Dm$.

Since we were careful to choose $E$ so that at least two points were
on each hemisphere, we can repeat this argument upside down to show
that $h_+$ is isotopic to power of $\TDp$ by an isotopy fixing the
boundary $\partial\Dp$ and all of the original $y_i$ on $\Dp$. By
joining these two isotopies of the hemispheres along $E = \partial\Dp
= \partial\Dm$, we acquire a isotopy on the sphere fixing all $n$
points.  Thus
$$
[h] = [(h_+,h_-)] = [(\TDp^r,\TDm^s)] = [T_E^{r+s}]
$$
in $\MSnp$, where $T_E$ is a Dehn twist along $E$, and $r$ and $s$ are
integers.

To finish the proof of the lemma, we note that $\Msfourp$ is a free
group; if we forget all but two points on each hemisphere, $[T_E]$ is
mapped to a generator of $\Msfourp$, and thus $r + s = 0$.  This is
where we need  that $n \geq 5$, so that there is point to be lost to 
obtain a trivial element of $\Msfourp$.  \endproof

We have shown that there are no nontrivial periodic or reducible elements
in $Br(S^2,n)$ for $n \geq 5$. Thus all nontrivial elements must be
pseudo-Anosov.

\end{section}

\section{Reducible and periodic mapping classes on the genus two surface} 
We now consider what happens to mapping classes of $\TT$ when they are
projected to the sphere.  Recall from section two that there is a
branched covering $\pi\co \TT \to S^2$ with branch points $y_1, \ldots,
y_6$, and that we have chosen these points to be the six specified
points for $\Ms$.  Also recall from Lemma 2.2 that there is a
surjective map from $\M$ to $\Ms$ whose kernel is generated by $[i]$,
the mapping class of the covering involution. The Brunnian subgroup
$\Brt$ is the preimage of $\Brus$ under this map.

Any periodic element of $\Brt$ projects to a periodic element of
$\Brus$, which must be trivial. Thus, the {\em only\/} periodic
element in $\Brt$ is the involution $[i]$.

\begin{thm}
  Reducible mapping classes in $\M$ project to reducible
  mapping classes in $\Ms$.
\end{thm}

Let $x_1,...,x_6 \in \TT$ be the preimages of the six branch points
 $y_1,..,y_6 \in S^2$.  Let $C$ be a simple closed curve that is
 invariant under the involution $i$ and contains $x_1$ and $x_2$.  $C$
 will project down to a simple line segment connecting $y_1$ and $y_2$
 on $S^2$.  Let $C'$ on $S^2$ be a simple closed curve that is the
 boundary of a small connected neighborhood of that line segment.  See
 Figure 2.

\begin{figure}[ht!]
\centerline{
\begin{tabular}{cc}
\small
%\ShowGrid 
\SetLabels 
\E(0.47*0.77){$S$}\\
\E(0.07*0.38){$C$}\\
\endSetLabels 
\AffixLabels{\includegraphics[height=1.5in]{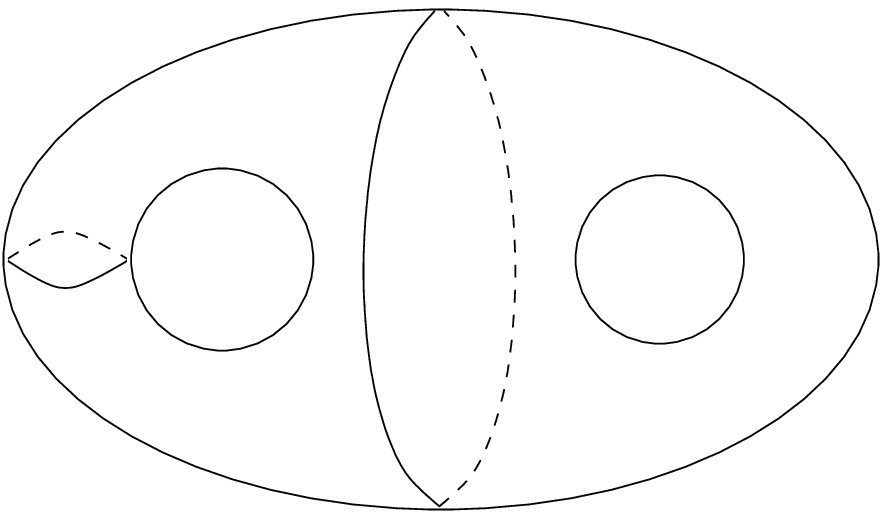}} 
& 
\small
%\ShowGrid 
\SetLabels 
\E(0.75*0.69){$S'$}\\
\E(0.21*0.70){$C'$}\\
\endSetLabels 
\AffixLabels{\includegraphics[height=1.5in]{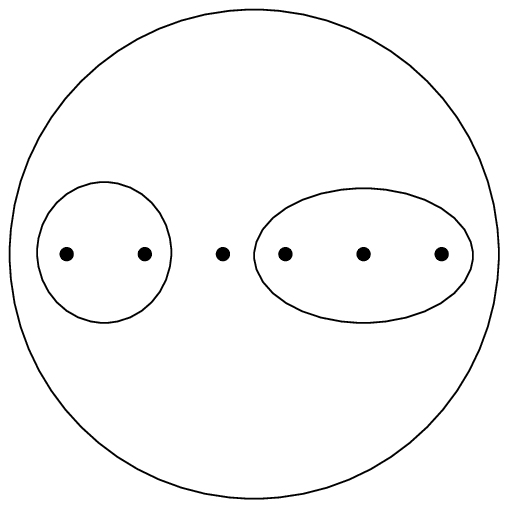}}
\\
(a) & (b)\\
\end{tabular}
}
\caption{(a) Example of $S$ and $C$\qquad\qquad  (b) Example of $S'$ and $C'$}
\end{figure}

Let $S$ be an essential separating simple closed curve that is also
 invariant under $i$, does not intersect $C$, and which separates the
 branch cover points into sets $\lbrace x_1,x_2,x_3 \rbrace$ and
 $\lbrace x_4,x_5,x_6 \rbrace$.  Let $S'$ on $S^2$ be the image of
 $S$.

Note that although there are many choices for $S$ and $C$, these
 choices differ only by the action of elements of $\M$.

\begin{lemma}
If $f\co \TT \to \TT$ is a homeomorphism that leaves $C$ invariant, then
$p([f])$ has a representative fixing $C'$.
\end{lemma}
\pf
The proof is to analyze what happens when we cut $\TT$ along $C$ and
then cap off one of the holes.

Let $N(C)$ be an annular neighborhood of $C$.  Then the closure of
$\TT \setminus N(C)$ is a genus one surface with two disjoint disks
removed.  We will call this surface $\Ttwo$. Recall that $\Mcc$ is the
group of isotopy classes of the set of homeomorphisms of $\Ttwo$ that
fix the boundaries pointwise.  Inclusion of $\Ttwo$ into $\TT$ induces
a well defined map
$$
 \phi\co   \Mcc \to \M.
$$
Given an orientation preserving homeomorphism $f$ that leaves $C$
invariant, we can assume that $f$ or $fi$ restricts to the identity on
the annular neighborhood $N(C)$.  Since we need only prove a fact
about $p([f])$, assume that $f$ fixes the annular neighborhood. Then
$f$ restricts to a homeomorphism {$f'\co \Ttwo \to \Ttwo$} that fixes the
boundary components pointwise.  Thus $[f']$ is an element of $\Mcc$,
and $\phi([f']) = [f]$; ie, $[f]$ is in the image of $\phi$.  To
prove the lemma it is enough to find generators for $\Mcc$ and show
that their images in $\M$ satisfy the conclusion.

Recall that we define $\Mccv$ to be the group of isotopy classes of the set of
homeomorphisms of $\Ttwo$ that may vary along the boundary components,
but do not switch them. 
 There is an exact sequence
$$
1 \to \tw \to \Mcc \to \Mccv \to 1
$$
where $\Tbp$ and $\Tbm$ are the twists about the boundaries.

If we cap off the holes of $\Ttwo$ with pointed disks, we obtain a
genus one torus $T$ with two specified points, $p_1$ and $p_2$.  Recall that
$\Mtp$ is the group of isotopy classes of the set of homeomorphisms of
$T$ that fix both $p_1$ and $p_2$.  Lemma 3.1 tells us that the
homomorphism from $\Mccv$ to $\Mtp$ induced by the capping off will be
an isomorphism.
$$
  \Mccv \cong \Mtp 
$$

If $L_2\co  \Mtp \to \Mt$ is the ``forget $p_2$'' map, then by the Remark
after Lemma 3.2, the following sequence is exact.
$$
 1 \to \pT \to \Mtp \stackrel{L_2}{\to} \Mt \to 1
$$
\begin{defn}\rm
The image of a standard generator of $\pT$ in $\Mtp$ is a
\textit{double Dehn twist,} that is, a pair of Dehn twists in opposite
directions on curves immediately to either side of some simple closed
curve representing the generator \cite{kn:jon}.
\end{defn}

Each of these double Dehn twists is isotopic to the identity by an
isotopy which fixes $p_1$ but lets $p_2$ move about; the tracks of
$p_2$ under these isotopies represent the corresponding elements of
$\pT$.

If we let $a$ and $b$ be the standard generators for $\pT$, then we
define $[TT_a]$ and $[TT_b]$ as the images of $a$ and $b$. 

$\Mt$ is generated by Dehn twists $[T_m]$ and $[T_l]$ along the
meridian and longitude of the punctured torus.  Thus $\Mtp$ is
generated by $[TT_a], [TT_b], [T_m]$, and $[T_l]$.  See Figure 3.  We
will denote the preimages of these isotopy classes in $\Mcc$ by
$[\widehat{TT}_a], [\widehat{TT}_b], [\widehat{T_m}]$, and
$[\widehat{T_l}]$.  These, and $[\Tbp]$ and $[\Tbm]$, generate $\Mcc$.

\begin{figure}[ht!]
\centerline{\includegraphics[width=\hsize]{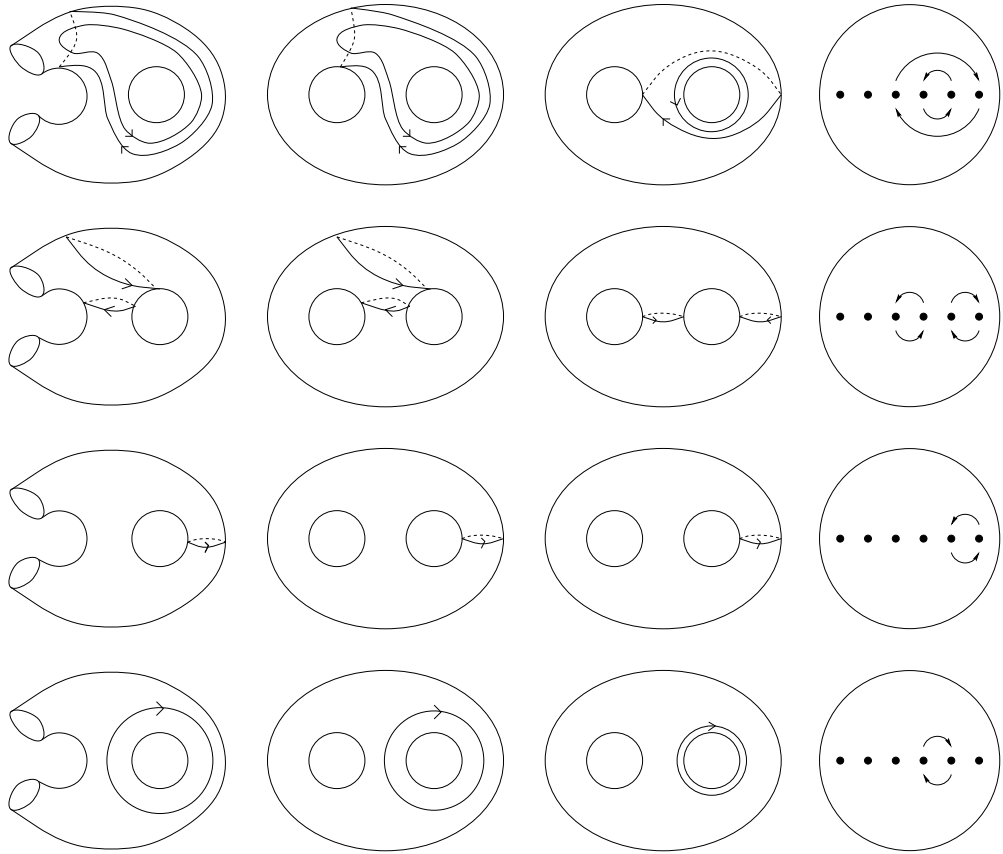}}
\caption{The four maps  $[TT_a], [TT_b], [T_m]$, and $[T_l]$, and their images in $\M$ and $\Ms$}
\end{figure}

We can easily see that the images of these mapping classes in $\M$ all
project down to mapping classes of the sphere with representatives
that fix $C'$.
\endproof

\begin{lemma}
If $f\co \TT \to \TT$ is a homeomorphism that leaves $S$ invariant, then
$p([f]^2)$ has a representative fixing $S'$.
\end{lemma}
\pf
Again, we cut along $S$ and analyze the result.  

Let $\torone$ and $\tortwo$ be the closures of the components of $\TT
\setminus S$. Then the mapping class groups $\Mtorone$ and $\Mtortwo$ are
generated by Dehn twists along their meridians and longitudes and by
Dehn twists along their boundaries.  Again, there is a well defined
map
$$
\phi \co  \Mtorone \times \Mtortwo \to \M,
$$
and the images of the generators project down to isotopy classes in
$\Ms$ that have representatives fixing $S'$. 

Given a homeomorphism $f$ that leaves $S$ invariant, $f$ either
switches the two components of $\TT - S$, or it does not.  In either
case, $f^2$ cannot switch the two components of $\TT - S$. Thus we can
assume that $f^2$ fixes an annular neighborhood of $S$.  Therefore
$f^2$ induces maps $f'$ and $f''$ on $\torone$ and $\tortwo$ which fix
the boundary, and hence $[f'] \in \Mtorone$ and $[f''] \in \Mtortwo$.
But then $[f]^2 = \phi([f']\times[f''])$.  \endproof

\pffive Let $[f] \in \M$ be a reducible mapping class. 
Then $[f]$ is conjugate to a mapping class $[g]$ which must fix one of
the following collections of nonintersecting simple closed curves:
\begin{enumerate}
\item {$S$}
\item {$S$ and $C$}
\item {$S, C$, and another essential non-separating curve}
\item {$C$}
\item {$C$ and one other essential non-separating curve}
\item {$C$ and two other essential non-separating curves}
\end{enumerate}
In cases 1, 2, and 3, $S$ must be left invariant, since it is the only
separating curve.  Therefore $p([g]^2)$ has a representative leaving
$S'$ invariant, by Lemma 5.3. In cases 4, 5, and 6, there is some
nonzero power $k$ so that $g^k$ has a representative leaving $C$
invariant.  But then $p([g]^k)$ leaves $C'$ invariant, by Lemma
5.2. Nonzero powers of reducible, periodic, and pseudo-Anosov elements
are reducible, periodic, and pseudo-Anosov, repsectively. Therefore,
by the classification theorem, in all cases above $p[g]$ is reducible,
and so is $p[f]$.

We have shown that reducible elements of $\M$ project to reducible
elements of $\Ms$. Therefore, since $\Brus$ has no nontrivial
reducible elements, $\Brt$ also has no nontrivial reducible elements
except the involution $[i]$.

\section{Nontriviality}

In Section 4 we proved that $\Br$ has no nontrivial periodic or
reducible elements.  Thus, all of its nontrivial elements are
pseudo-Anosov. In Section 5, we proved that the only nontrivial
element of $\Brt$ that is not pseudo-Anosov is the involution
$[i]$. To rid ourselves of this element, we intersect with the kernel
of the usual $\Z_3$ homology representation, $\rho\co  \M \to \slthree$.

\begin{thm}{\rm(Long \cite[page 83]{kn:lon})}\qua
The intersection of two nontrivial, noncentral normal subgroups of a
mapping class group of a surface is nontrivial.  
\end{thm}
Since $\Brt$ and $\ker{\rho}$ are not central, their intersection in
$\M$ is nontrivial.  As an example of a nontrivial element of $\Brt
\cap \ker{\rho}$, consider the nested commutator
$$
  [\Delta_1^6,[\Delta_2^6,[\Delta_3^6,[\Delta_4^6,\Delta_5^6]]]],
$$
where $\Delta_i$ are the standard generating Dehn twists on $\TT$.
This mapping class is Brunnian, since it projects to
$$
 [\sigma_1^6,[\sigma_2^6,[\sigma_3^6,[\sigma_4^6,\sigma_5^6]]]],
$$
which clearly reduces to the identity if any of the six points on the
sphere is forgotten.  It also acts trivially on homology with $\Z_3$
coefficients, but it is
\textit{not} null-homologous over $\Z$, and so it is not trivial.  

Other examples can be easily generated by taking commutators of
Brunnian mapping classes and null-$\Z_3$--homologous elements, or
simply by using homology to check that a given Brunnian mapping class
is not $[i]$.  This provides a new source of pseudo-Anosov mapping
classes of the genus two torus.

\rk{Acknowledgments}
This paper is part of my dissertation. I would like to thank John
Stallings for serving as my thesis advisor, and Andrew Casson for
numerous conversations and suggestions.  I would also like to thank
the referees. This research was partly funded by NSF grant DMS
950-3034.


\begin{thebibliography}

\bibitem{kn:birm}{\bf J Birman},  
\emph{Braids, Links, and Mapping Class groups},
Annals of Math. Studies, No. 82. Princeton University Press

\bibitem{kn:birm2}{\bf J Birman}, {\it Mapping class groups and their
relationship to braid groups},  Comm. Pure and Applied Math. 22
(1969) 213--238 

\bibitem{kn:edkr}{\bf R Edwards}, {\bf R Kirby}, {\it Deformations of spaces
of imbeddings},  Annals of Math. 98 (1971) 63--88 

\bibitem{kn:flp}{\bf Fathi}, {\bf F Laudenbach}, {\bf V Ponerau}, et al,
\emph{Travaux de Thurston sur les surfaces},
Ast\'erisque 66--67 (1979)

\bibitem{kn:jon}{\bf D Johnson}, {\it Homeomorphsims of a surface
which act trivially on homology}, Proc. of AMS, 75 (1979) 119--125

\bibitem{kn:lev}{\bf H Levinson}, {\it Decomposable braids and linkages}, 
Trans. AMS, 178 (1973) 111--126  

\bibitem{kn:lon}{\bf D Long}, {\it A note on the normal subgroups of 
mapping class groups} Math. Proc. Camb. Phil. Soc. 99 (1985) 79--87

\bibitem{kn:rolf}{\bf D Rolfsen}, \emph{Knots and Links}, Publish or
Perish, Inc.  Math. Lect. Series 7 (1976)

\end{thebibliography}
\end{document}